# Sieve estimates for biased survival data[*]

## Jiayang Sun[1] and Bin Wang[2]

*Case Western Reserve University and University of South Alabama*

**Abstract:** In studies involving lifetimes, observed survival times are frequently censored *and* possibly subject to biased sampling. In this paper, we model survival times under biased sampling (a.k.a., biased survival data) by a semi-parametric model, in which the selection function $w(t)$ (that leads to the biased sampling) is specified up to an unknown finite dimensional parameter $\theta$, while the density function $f(t)$ of the survival times is assumed only to be smooth. Under this model, two estimators are derived to estimate the density function $f$, and a pseudo maximum likelihood estimation procedure is developed to estimate $\theta$. The identifiability of the estimation problem is discussed and the performance of the new estimators is illustrated via both simulation studies and a real data application.

## 1. Introduction

The problem of analyzing survival data arises in many application fields, such as clinical trials in medicine, reliability assessments in engineering, biology, epidemiology and public health. Censoring is a common phenomenon accompanying survival data, often due to voluntary or involuntary drop-out of study subjects. In addition, survival data may have been drawn by *biased sampling* (with or without our knowledge) in which whether a survival time $T$ can be observed depends on a *selection function* $w(t)$ which is the probability of observing $T$ if the true value of $T$ is $t$. Survival data drawn under such biased sampling, when $w(t)$ is not a constant, are hereafter called *biased survival data/sample*. When $w(t)$ is a constant, the survival data are called the *standard survival data*. Here are three examples of biased survival data, with the first given in more detail.

1. In a study of Scleroderma, a rare disease, some data of all cases of Scleroderma diagnosed in Michigan from 1980 to 1991 were collected and the times from diagnosis of Scleroderma to death were recorded [2]. Based on the Kaplan–Meier (K–M) estimates of survival curves for patients diagnosed in 1980–1985 versus 1986–1991, Gillespie (one of the authors in [2]) found that the earlier group patients (from 1980–1985) lived significantly longer than the later group patients (from 1986–1991). What had happened? If anything had changed, medical care should have improved in 1986–1991 over 1980–1985 and hence the second group of patients should have had better survival times. According to Gillespie, their sources of information included hospital databases and responses from private physicians. Unfortunately, because hospital records did not always go back to 1980, and physicians did not always remember patients they saw many years ago, patients who were still alive (and thus had more current hospital records) were more likely to be collected in the sample than those who died in the early period. This resulted in a biased survival sample for the 1980–1985 group. Indeed, as Gillespie stated "We feel that the result is entirely due to our length-biased sample." Length-biased sampling is a special example of biased sampling with $w(t) \propto t$.

---

[*]The Research is supported in part by an NSF award.

[1]Case Western Reserve University, e-mail: jiayang@sun.cwru.edu
[2]University of South Alabama, e-mail: bwang@jaguar1.usouthal.edu

2. In assessing familial risk of disease based on a reference database, which is a collection of family histories of cases typically assembled as a result of one family member being diagnosed with a disease. Clearly, the larger a family is, the greater the probability that this family will be found from the registry is [3].

3. In cancer screening programs, whether the pathological changes of a patient in the preclinical phase can be discovered depends very much on the phase of the tumor.

Let $f(t)$ be the true probability density function (pdf) of the survival time $T$ and $F(t)$ the corresponding cumulative distribution function (cdf). If the sampling bias in a biased survival sample is ignored in an estimation of $f$ or $F$, resulting estimates are not consistent and can be misleading as shown in Example 1. In fact, the missingness resulting from biased sampling is sometimes also called "non-ignorable" missing because it leads to an observed sample that has a density weighted by $w$, as demonstrated in (1) below, in contrast to missing at random, MAR, in which whether a subject is missing is independent of $t$ and may be ignored.

In this paper, we propose a semi-parametric model that incorporates both the censoring information and biased sampling scheme for modeling biased survival data (Section 2). In our model, the density function $f(t)$ of the survival times is assumed to be smooth, and the selection function $w(t)$ is specified up to an unknown finite dimensional parameter $\theta$ and is a constant when $\theta = \theta_0$; for example, $w(t) \propto t^\theta$, for $\theta \geq 0$. So this model is applicable to both biased survival data and standard survival data. The identifiability of estimating $(f, w)$ is also discussed. The semi-parametric parameter $(f, \theta)$ under our semi-parametric model is "sieve" identifiable. In Section 3, two estimators, one called weighted kernel estimator (WKE) and the other called transformation-based estimator (TBE), are derived for estimating $f$; and a "pseudo" maximum likelihood procedure is proposed for estimating $\theta$. Our new estimators are compared with those that ignore either censoring or sampling biases or assume that the selection function is known; and examined as the sample size increases (Section 4). The $L_1$, $L_2$ distances and MSE of our estimators $\hat{f}$ converge, while the naive estimator (that ignore both censoring and selection biases), the K–M estimator and the Jones estimator [8] did not perform as well as ours. In terms of a confidence interval for $F(t)$, our WKE and TBE also beat the naive, the K–M and Jones estimators. The application of the new estimators is illustrated via an analysis of a survival data set on time until death of bone-marrow transplant patients in Section 5. The paper concludes with some discussions in Section 6.

## 2. Model

In an idealized situation, there would be $T_1, \ldots, T_N \stackrel{iid}{\sim} f$ for us to make an inference about $f$. In reality, we observe only a subset of $T_1, \ldots, T_N$, where each $T_i$ is included in the subset with probability $w(t_i)$ if the value of $T_i$ is $t_i$. The function $w(t)$ is called the *selection function*. Abusing the notation a little, we denote the subset by $\mathcal{T} = \{t_1, t_2, \ldots, t_n\}$, then the observed sample size $n \sim$ Binomial $(N, \kappa)$, where $\kappa = E_f(w(T)) < \infty$ is the mean value of the probability of observing $T$. The observed sample $\mathcal{T}$ no longer has the common pdf $f$. Instead, conditioning on $n$,

$$(1) \qquad t_1, \ldots, t_n \overset{iid}{\sim} f_w(t) \equiv \frac{w(t)f(t)}{\kappa}.$$

Thus, if a standard procedure that ignores the selection bias in $\mathcal{T}$ is used, the resulting density estimate based on $t_i$'s might be consistent to $f_w$, but not to $f$, the density of interest. A new procedure that accounts for both censoring and selection biases must be developed to delineate $f$ and $w$.

*Identifiability.* If both $f(t)$ and $w(t)$ in (1) are completely unknown, the problem of estimating both $f$ and $w$ based on $\mathcal{T}$, *one* biased sample, is unidentifiable. For any $w(t)$ and an arbitrary $h(t) > 0$ for which the integral $\kappa' = \int f(u)/h(u)du$ is finite, the pair $(w(t), f(t))$ gives the same likelihood as the pair $(h(t)w(t), f(t)/(h(t) \cdot \kappa'))$. To ensure the identifiability, (a) either $w(t)$ or $f(t)$ has to be assumed to be parametric; or (b) there is another biased sample $\mathcal{T}'$ which has some overlap with $\mathcal{T}$. The overlap $\mathcal{T} \cap \mathcal{T}'$ provides some information about $w$ and hence allows both $f$ and $w$ to be estimable nonparametrically in the range of $\mathcal{T} \cap \mathcal{T}'$. See Lloyd and Jones [12], Wang and Sun [25], for more information on nonparametric estimates of $f$ and $w$, based on *two* biased samples.

In this paper, we consider the case when there is only *one* biased survival sample, which is $\mathcal{T}$ with some of its members censored. Here $N$, the size of the idealized sample is not assumed to be known. In other words, our first model assumption is:

**Model Assumption 2.1.** The observable sample is: $(T, I) = \{(\tilde{t}_i, I_i), i = 1, \ldots, n\}$, where $\tilde{t}_i = t_i \sim f_w$ if $t_i$ is uncensored $(I_i = 1)$, and $\tilde{t}_i = c_i$, a censoring time, if $t_i$ is censored $(I_i = 0)$. The censoring times $c_i$'s are independent of survival times $t_i$'s and have a common censoring distribution that has the same support as that of $f_w$. Further, a right censoring scheme is assumed: $I_i = 1$ if $t_i \leq c_i$, 0 otherwise.

To ease the notation, we abuse the notation again to simplify $\tilde{t}_i$ by $t_i$ hereafter. So, we observe $\{(t_i, I_i)\}$, where if $I_i = 1$, the actual survival time is $t_i \sim f_w$, uncensored, if $I_i = 0$, the survival time is censored and is greater than $t_i$.

Next, given only *one* biased survival sample, we shall assume either $w$ or $f$ to be a parametric function. If both $f$ and $w$ are parametric functions, estimating $f$ and $w$ is equivalent to estimating parameters in a parametric survival model. For example, let

$$f(t) \propto \gamma \alpha t^{\alpha-1} e^{-\gamma t^{\alpha}}, \quad \text{and } w(t) \propto t^{\beta},$$

then the weighted density is

$$f_w(t) = \frac{\gamma \alpha}{\kappa'} t^{\alpha+\beta-1} e^{-\gamma t^{\alpha}}.$$

The unknown parameters, $\alpha, \beta, \gamma$, can be estimated by maximizing the likelihood:

$$(2) \qquad L = \prod_i [f_w(t_i)]^{I_i} [S_w(t_i)]^{1-I_i},$$

where $S_w(t_i) = 1 - \int_0^{t_i} f_w(u)du$ is the survival function at $t_i$. The expression of the resulting mle from (2) may be complicated but the estimates are straightforward to compute. Hence, as long as the parameters are identifiable from (2), the parameters can be estimated using standard parametric estimation procedures.[1]

---

[1]In the case where both $f(t)$ and $w(t)$ are specified up to an unknown finite dimensional parameter, typical identifiability conditions for a parametric model are still needed (though standard) to estimate the parameters consistently. For example, if $f(t) \propto e^{\alpha t}, w(t) \propto e^{-\beta t}$, then $f_w(t) \propto e^{(\alpha-\beta)t}$. Clearly, only $\alpha - \beta$ is identifiable based on one biased sample.

The estimation problem becomes more interesting and challenging when $f(t)$ is assumed to be smooth while $w(t)$ is a parametric function. This semi-parametric model is more general than a parametric model and is useful when there is no obvious choice for a parametric assumption on $f$. So, we assume next:

**Model Assumption 2.2.** For the biased survival data in Model Assumption 2.1, the pdf $f$ is assumed to be smooth and the selection function $w(t)$, denoted hereafter by $w(t, \theta)$, is specified up to an unknown finite dimensional parameter $\theta \in \Theta$. Hence the weighted density is now

$$(3) \qquad\qquad f_w(t) = \frac{w(t, \theta) f(t)}{\kappa(\theta)}$$

where $\kappa(\theta) = E_f(w(T, \theta))$ with the expectation taken for $T \sim f$.

The semi-parametric model specified by (3) is related to those considered by Gill et al. [6], Robbins and Zhang [14], and Vardi and Zhang [21], among others. A notable difference is that in our semi-parametric model (which satisfies Model Assumptions 2.1 and 2.2), $f$ is assumed smooth so that nonparametric smoothing techniques can be used in estimating $f$ and a sieve estimate of $\theta$ based on a pseudo likelihood can be developed as that in Section 3.

Let $S(\theta) = \{t : w(t, \theta) > 0\}$ be the support of $w(t, \theta)$. If $S(\theta)$ depends on $\theta$ in that both sets $S(\theta) - S(\theta')$ and $S(\theta') - S(\theta)$ have positive measure under $F$ for all $\theta \neq \theta' \in \Theta$, then both $\theta$ and $f$ are completely identifiable as shown by Lemma 2.3 of Gilbert et al. [5]. In practice, the selection function may be a polynomial function in $t$, e.g. $w(t, \theta) \propto t^\theta, \theta \geq 0$, and $w(t, 0)$ is a constant. This $w$ has support $(0, \infty)$, which is independent of $\theta$. So it does not satisfy the condition of this Lemma 2.3. However, we can put a constraint on the form of $w(t, \theta)$ and the type of $\widehat{f}(t)$, for each fixed sample size. Then the resulting semi-parametric estimator, under a sieve identifiability defined in §3, will be similar to those obtained by "sieve" methods and hence will lead to reasonable estimators of $\theta$ and $f$. See Section 3.

Alternatively, if one can model $f$ as a parametric function, the assumption for $w$ can be relaxed to be nonparametric. In Sun and Woodroofe [18], $f$ is assumed to come from an exponential family of distributions and $w$ is assumed only to be monotone. They also developed an iterative MM (maximization-maximization) algorithm for estimating both $w$ and parameter in $f$ when $N$ in the idealized situation is known and unknown (two very different cases). They showed that the MM algorithm converges to correct (penalized) maximum likelihood estimators and the estimators are consistent. This type of semi-parametric model is dual to the semi-parametric model proposed above and may be extended to allow for censored observations. We do not consider this extension in this paper. For a recent tutorial on MM algorithms under other settings, see [7].

## 3. Semi-parametric estimators

In this section, we develop semi-parametric estimators of $(f, \theta)$ under Model Assumptions 2.1 and 2.2, and discuss an additional identifiability condition required for our estimation procedure.

### 3.1. Weighted kernel estimator (WKE)

The bias due to censoring can be corrected in a standard kernel estimator by weighting the K–M estimator with the censoring information as proposed by Marron and Padgett [13]. The basic idea is as follows: Order the sample $(T, I)$ with respect to $T$ and denote it by $\{(t_{(i)}, I_{[i]}), i = 1, \ldots, n\}$. Then the K–M estimate of the cdf is

$$(4) \qquad \widehat{F}_{km}(t) = \begin{cases} 0, & 0 \le t \le t_{(1)}, \\ 1 - \prod_{i=1}^{j-1} \left( \frac{n-i}{n-i+1} \right)^{I_{[i]}}, & t_{(j-1)} < t \le t_{(j)}, \\ 1, & t > t_{(n)}. \end{cases}$$

The Marron and Padgett's kernel estimator of $f(t)$ induced by $\widehat{F}_{km}$ is then

$$\widehat{f}(t; h) = \int K_h(t - z) d\widehat{F}_{km}(z) = \sum_{i=1}^{n} s_i K_h(t - t_{(i)}),$$

where $K_h(t) = (1/h)K(t/h)$, $K$ is a symmetric probability density kernel such as the $N(0, 1)$ density function, and $s_i$ is the size of the jump of $\widehat{F}_{km}$ in (4) at $t_{(i)}$.

We can correct the selection bias by replacing the weight function $s_i$ with $s_i / w(t_{(i)}, \theta)$. Therefore, a new *weighted kernel estimator* is proposed,

$$(5) \qquad \widehat{f}_{wk}(t) = \widehat{\kappa}_{wk} \sum_{i=1}^{n} \frac{s_i K_h(t - t_{(i)})}{w(t_{(i)}, \hat{\theta})},$$

where $\widehat{\kappa}_{wk}$ is a normalizing constant, such that $\widehat{\kappa}_{wk}^{-1} = (\sum s_i / w(t_{(i)}, \hat{\theta}))$, and $\hat{\theta}$ is a good estimate of $\theta$, such as one described in Section 3.3. If $w(t_{(i)}, \hat{\theta}) \propto t_{(i)}$ is a known length-biased selection function, $\widehat{f}_{wk}(t)$ in (5) is reduced to the Jones estimate [8]. See the comparisons of $\hat{f}_{wk}$ with the Jones estimate in Section 4.

### 3.2. Transformation based estimator (TBE)

Another way to correct both the selection and censoring biases is by using the transformation-based method of El Barmi and Simonoff [4] to correct for the selection bias, and using $s_i$ from the K–M estimate of the transformed variable to account for the censoring bias, simultaneously.

Let $g(y)$ be the density function of $Y = W(T) \equiv W(T, \theta)$, where $T \sim f_w$ and

$$W(t, \theta) = \int_0^t w(u, \theta) du$$

is the *cumulative selection function*. For example, if $w(t, \theta) = c \cdot t^{\theta}$, for a constant $c > 0$ and $\theta \ge 0$, then $W(t, \theta) = c \cdot t^{\theta+1}/(\theta + 1)$ is monotone in $t$ on $[0, \infty)$. The cumulative distribution function of $Y$ can be easily shown to be

$$(6) \qquad G(y) = F_w\left(W^{-1}(y)\right),$$

where $W^{-1}(t)$ is the inverse function of $W(t, \theta)$ for fixed $\theta$ and $F_w(x) = \int_0^x f_w(u) du$ is the cdf of $f_w$. Differentiating $G(y)$, we obtain the pdf $g(y) = f(W^{-1}(y))/\kappa$. Thus,

$$f(t) = \kappa \cdot g(W(t, \theta)).$$

Hence, for fixed $\theta$, let $(Y, I) = \{(Y_i, I_i), i = 1, \ldots, n\}$, where $Y_i = W(t_i, \theta)$. Order this sample with respect to $Y$ and denote it by $\{(Y_{(i)}, I_{[i]}))\}$. Then the K–M estimator of the cdf of $Y$ is

$$\widehat{F}_{km}(y) = \begin{cases} 0, & 0 \le y \le Y_{(1)}, \\ 1 - \prod_{i=1}^{j-1} \left( \frac{n-i}{n-i+1} \right)^{I_{[i]}}, & Y_{(j-1)} < y \le Y_{(j)}, \\ 1, & y > Y_{(n)}. \end{cases}$$

Let $s_i$ denote the jump size of this K–M estimate at $Y_{(i)}$. Then our proposed *transformation based estimator* is

$$(7) \qquad \widehat{f}_{tb}(t) = \widehat{\kappa}_{tb} \sum_{i=1}^{n} s_i K_h(W(t, \theta) - Y_{(i)}),$$

where $\widehat{\kappa}_{tb}$ is a normalizing constant such that $\widehat{\kappa}_{tb}^{-1} = \sum s_i K_h(W(t, \theta) - Y_{(i)})$. Here $\theta$ is replaced by a good estimate $\widehat{\theta}$ when $\theta$ is unknown. See next section for an estimate of $\theta$. If $\theta$ is known and $s_i = 1/n$ for all $i$, $\widehat{f}_{tb}(t)$ is reduced to the El Barmi and Simonoff estimate.

### 3.3. Estimator of θ

If $\theta$ is unknown, we propose to estimate it by maximizing a corresponding "pseudo" or "sieve" log-likelihood:

$$(8) \qquad \widehat{l}_{wk}(\theta) = \sum_{j=1}^{n} I_j \log[\widehat{f}_{wk}(t_j, \theta)] + \sum_{j=1}^{n} (1 - I_j) \log[\widehat{S}_{wk}(t_j, \theta)], \quad \text{or}$$

$$(9) \qquad \widehat{l}_{tb}(\theta) = \sum_{j=1}^{n} I_j \log[\widehat{f}_{tb}(t_j, \theta)] + \sum_{j=1}^{n} (1 - I_j) \log[\widehat{S}_{tb}(t_j, \theta)]$$

where $\widehat{f}_{wk}(t, \theta) = \widehat{f}_{wk}(t)$ and $f_{tb}(t, \theta) = \widehat{f}_{tb}(t)$ are defined in (5) and (7) with $\widehat{\theta}$ replaced by $\theta$, and $\widehat{S}_{wk}(t_j, \theta)$ and $\widehat{S}_{tb}(t_j, \theta)$ are the survival functions at point $t_j$ for the two methods respectively,

$$\widehat{S}_{wk}(t, \theta) = 1 - \int_0^t \widehat{f}_{wk}(u, \theta) du, \qquad \widehat{S}_{tb}(t, \theta) = 1 - \int_0^t \widehat{f}_{tb}(u, \theta) du.$$

In the rest part of this paper, the following "sieve" identifiability is assumed:

**Model Assumption 3.1 (Sieve identifiability).** The semi-parametric model with unknown parameters $\theta$ and $f$ is "sieve" identifiable in the following sense:

$$\widehat{l}_{wk}(\theta_1) = \widehat{l}_{wk}(\theta_2) \quad \text{for a.s. all } t_i \in \mathcal{R}^+ \quad \Longleftrightarrow \quad \theta_1 = \theta_2,$$
$$\widehat{l}_{tb}(\theta_1) = \widehat{l}_{tb}(\theta_2) \quad \text{for a.s. all } t_i \in \mathcal{R}^+ \quad \Longleftrightarrow \quad \theta_1 = \theta_2,$$

where $\mathcal{R}^+$ is the support of $f$. For practical purposes, in the one-dimensional case, $R^+$ can be taken as $(0, a)$ for some large $a > 0$.

This type of identifiability ensures that $\theta$ is identifiable under the sieve likelihood (8), or (9), respectively, and the mle of $\theta$ from the corresponding sieve likelihood exists. Call the $\widehat{\theta}$ which maximizes (9), the *pseudo mle*. Since the sieve likelihood is usually a good approximation to the true likelihood as $n \to \infty$, we expect our WKE and TBE $\widehat{f}$ based on the $\widehat{\theta}$ to be consistent. This is very much in the same spirit as that of a histogram estimator. A properly chosen histogram estimator is consistent to $f$ under some regularity conditions while the fully nonparametric mle of $f$ is a useless and inconsistent estimate of $f$. The fully nonparametric mle places a delta function at every data point. The consistency of our WKE and TBE are confirmed by Table 1 in Section 4. Our final WK and TB estimators of $f$ are $\widehat{f}_{wk}(t, \widehat{\theta}_{wk})$ and $\widehat{f}_{tb}(t, \widehat{\theta}_{tb})$, where $\widehat{\theta}_{wk}$ and $\widehat{\theta}_{tb}$ are the respective pseudo mle's from the corresponding

WK and TB sieve likelihoods (8) and (9). See Section 6 for a further discussion on the asymptotic justification of our proposed procedures.

In some extreme cases or when the sample size is not large enough, the optimal value of $\theta$ may be located at the edge of the specified range $\Theta$ of $\theta$. The penalized log-likelihoods of the form

$$(10) \qquad \widehat{\ell}^*_{wk} = \log \widehat{L}_{wk} = \sum_{j \in \mathcal{U}} \log \widehat{f}_{wk}(t_j, \theta) + \sum_{j \in \mathcal{C}} \log \widehat{S}_{wk}(t_j, \theta) - \frac{\alpha n}{\kappa_{wk}},$$

$$(11) \qquad \widehat{\ell}^*_{tb} = \log \widehat{L}_{tb} = \sum_{j \in \mathcal{U}} \log \widehat{f}_{tb}(t_j, \theta) + \sum_{j \in \mathcal{C}} \log \widehat{S}_{tb}(t_j, \theta) - \frac{\alpha n}{\kappa_{tb}},$$

are then considered to overcome this difficulty, where $\mathcal{C} = \{i : I_i = 0\}$ and $\mathcal{U} = \{i : I_i = 1\}$, $0 < \alpha \leq 1$ may approach zero as $n \to \infty$, which was discussed in details by Woodroofe and Sun [26]. This penalized log-likelihood is maximized subject to the constraint

$$w(t, \theta) \geq \epsilon, \quad \text{for all } t \in \mathbf{R}^+, \epsilon \geq 0 \text{ and } \sup_t w(t, \theta) = 1.$$

Under this constraint, "$w(t, \theta) \propto t$" means that $w(t, \theta)$ is only proportional to $t$ in $(\epsilon, a - \epsilon) \in R^+$ for some $\epsilon > 0$.

In this study, we take $\alpha = cn^{-0.5}$, where $c$ is a constant and its value can be chosen by the Jackknife or Cross-validation method. In this paper, we choose $c$ by minimizing either of the following expressions,

$$(12) \qquad CV_1 = \frac{n-1}{n} \sum (\widehat{\theta}_{-i,c} - \widehat{\theta}_{.,c})^2 + (n-1)^2 (\widehat{\theta}_c - \widehat{\theta}_{.,c})^2,$$

$$(13) \qquad CV_2 = \frac{1}{n} \sum_i \widehat{f}_{-i,c}(t_i),$$

$$(14) \qquad CV_3 = \frac{1}{\widehat{\kappa}_c (1 - \widehat{\kappa}_c)} \left\{ \frac{n-1}{n} \sum (\widehat{\kappa}_{-i,c} - \widehat{\kappa}_{.,c})^2 + (n-1)^2 (\widehat{\kappa}_c - \widehat{\kappa}_{.,c})^2 \right\},$$

where the subscript "$-i$" means that the $i$th data point has been omitted and the subscript "$.$" denotes an average of $*_{-i}$'s. The CV estimation of $c$ can be computationally intensive. For large data sets, the fast Fourier transformation may be implemented to speed up the algorithm [17].

## 4. Simulation Studies

### 4.1. Setup

In this simulation study, we consider a Weibull density with shape parameter $\gamma = 2$ and scale parameter $\lambda = 1$,

$$(15) \qquad f(t, \gamma, \lambda) = \lambda \gamma t^{\gamma - 1} \exp(-\lambda t^\gamma).$$

The solid line in Figure 1 shows the density curve of $f$ defined in (15).

To show the results of ignoring either sampling or censoring biases in a typical density estimate, we draw four samples using the following four designs. The kernel density curves of these four samples (without a correction for either selection or censoring biases), $\widehat{f}(t) = 1/(nh) \sum_i K((t - t_i)/h)$ for $t_i \in S$, are shown in Figure 1.

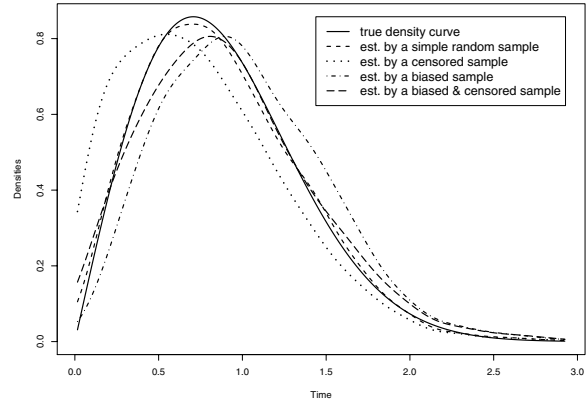

FIG 1. *Kernel estimates of densities for samples with sampling bias and/or censoring.*

- Simple random sample ($S = \mathbf{S}_0$): A simple random sample of size 3000 was drawn from $f$. The density curve of this sample was estimated by using the standard kernel method and is shown by the short-dashed curve in Figure 1. It's easy to see this curve is close to the true density curve of $f$ as expected.
- Sample with censoring bias only ($S = \mathbf{S}_c$): A sample of size 3000 was randomly selected from $f$ and 30% of the data points were randomly censored. As shown in Figure 1, the kernel density curve (dotted curve) of sample $\mathbf{S}_c$ shifts to the left of the true density curve $f$. This is also expected as now the sampling distribution (of $S_c$) is different from the target distribution ($f$). This is typical for right-censored survival data.
- Sample with selection bias only ($S = \mathbf{S}_b$): A sample of size 3000 was randomly chosen from the population. Each of these 3000 elements was observed subject to the selection probability $w(t, \theta) = w(t, 0.5) \propto \sqrt{t}$. This $w(t, \theta)$ implies that the elements with longer survival times were more likely to be sampled. The kernel density estimate of the density curve of sample $\mathbf{S}_b$ was computed and is shown as the dash-dotted line in Figure 1. We see that the sample density curve shifts to the right. This is also a case that the sampling distribution is different from the target distribution.
- Sample with both censoring and selection biases ($\mathbf{S}_{cb}$): In sample $\mathbf{S}_b$, if 30% of the data are further randomly censored, we obtain a biased survival sample. The density curve of $\mathbf{S}_{cb}$ is estimated and shown as the long-dashed curve in Figure 1. We find that the selection bias in this case somehow has balanced out the left-shift-ness of the density curve of $S_c$ though it is still not as good as the estimate based on a simple random sample from $f$. We can not rely on this kind of cancellation. If $w(t, \theta)$ had decreased with the increase of $t$, the selection bias would make the sample density curve more right-skewed.

The observed sample sizes $n$ were governed by the selection function and censoring scheme; they varied from one realization to another.

The results from these four experiments show that if a sampling distribution is different from a target distribution, then the deviation of the sampling distribution from the target distribution must be considered in developing a good estimate of the target density function, otherwise the resulting estimator is inconsistent.

### 4.2. Estimates based on a biased and censored sample

Using a biased sample that has some of the data points right-censored, we can estimate $f$ and $w$ by *WKE* and *TBE*. First, we estimate the unknown parameter $\theta$ in the selection function $w(t, \theta)$ by maximizing the log-likelihoods in (8) and (9) or the penalized log-likelihoods in (10) and (11).

Figure 2 shows the pseudo maximum likelihood estimates of the unknown parameter $\theta$. By using the WKE, we obtained an estimate of $\widehat{\theta} = 0.42$ (plot 2(a)), which is closer to the true $\theta = 0.5$ than that by using the TBE (plot 2(b)). We can then estimate $f$ by using the estimates in (5) and (7), by replacing $\theta$'s with $\widehat{\theta}$'s. In Figure 3, the thick solid line shows the true density curve of $f$, while the dashed line shows the WKE by treating the true parameter $\theta$ as known ($\theta = 0.5$) and the thin solid line shows the WKE by using the estimated parameter $\widehat{\theta} = 0.42$. The kernel density curve of the sample is also plotted as the dot-dashed line. From this figure, we can see that the three density curves are close and the WKE's based on known $\theta$ and $\widehat{\theta}$ are only slightly better. However, this result is based on only one

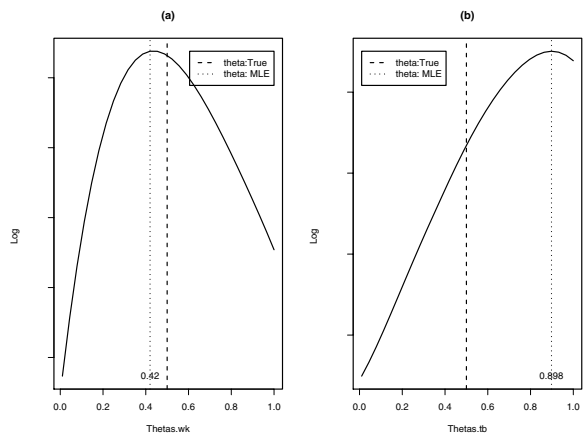

Fɪɢ 2. *Maximum likelihood estimators of* $\theta$.

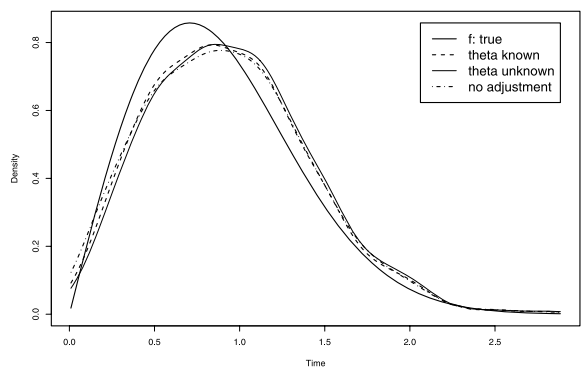

Fɪɢ 3. *Weighted kernel estimates of* $f$.

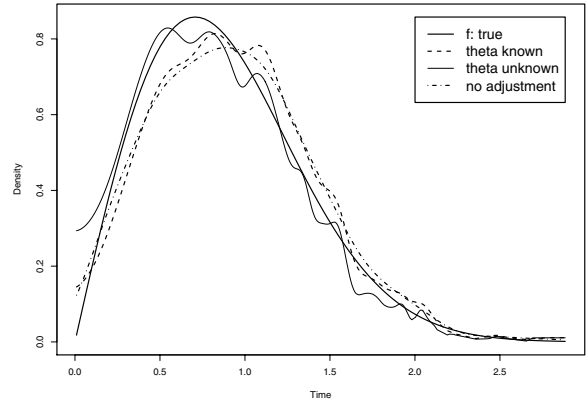

FIG 4. *Transformation based estimates of f.*

sample. See the next subsection for a report on the overall performance.

The TBE's obtained by using the true and estimated $\theta$ are displayed in Figure 4. From Figure 4, we can see that by using the true parameter $\theta$, we obtained an estimate (dashed line) which is close to the density curve of the sample (biased with censoring). While by using the estimated parameter $\widehat{\theta} = 0.898$, we obtained an estimate which is much closer to the true density curve $f$. The reason for this superiority is perhaps that fixing $\theta = 0.5$ may have limited some degrees of freedom of the semi-parametric approach – data can speak for themselves. The curves of TBE's are closer to the true density curve but coarser than those of WKE's. This is expected because the TBE corrects for the selection bias and censoring bias exactly in the same order as how biased survival data are formed, and the coarseness may come from the transformation or the way the window-width is determined. Further improvement is possible by applying some smoothing techniques to the TBE.

Hence, as an estimate of $f$, TBE estimate is the winner though it is a bit rough, but it can be smoothed out one more time.

## 4.3. Overall performance of WKE and TBE

To study the overall performance of the weighted kernel estimator and transformation-based estimator, we designed two experiments. The first experiment has the following design:

**Step 1:** Draw a sample $S$ of size $N = 50$ from $f$, subject to biased sampling with a selection function $w(t, \theta) = w(t, 1.0)$, and with 30% of the data points censored.

**Step 2:** Based on this sample $S$, estimate the cdf $F(t)$ by using the WKE, TBE, Jones estimate, naive estimate (with which we estimate the density function from the biased survival data, but without considering either the selection or censoring biases), denote the results by $\widehat{F}_{wke}$, $\widehat{F}_{tbe}$, $\widehat{F}_{jones}$, $\widehat{F}_{naive}$ respectively.

**Step 3:** Repeat step 1 and step 2 for 1000 times. Compute $L_1$, $L_2$ and MSE as

defined by

$$L_1 = \sum_i |\widehat{F}(t_i) - F(t_i)| \cdot d_i, \quad L_2 = \sum_i (\widehat{F}(t_i) - F(t_i))^2 \cdot d_i,$$

$$MSE = \frac{1}{n} \sum_i (\widehat{F}(t_i) - F(t_i))^2,$$

where $d_i = (t_{i+1} - t_{i-1})/2$ for $i = 2, 3, \ldots, n-1$ and $d_1 = t_2 - t_1$, $d_n = t_n - t_{n-1}$. So, $L_1$ and $L_2$ above are approximations of the $L_1$ and $L_2$ distances (in the form of integrals).

**Step 4:** Take other values of $N$, $N = 100, 200, \ldots$, and repeat step 1 through step 3.

Table 1 shows the $L_1$, $L_2$ distances, and MSE of $\hat{f}$ from the true $f$. Note that since the real $\theta = 1$ in this case, the assumptions used in the Jones estimate are justified and the Jones estimate is equivalent to the WKE with a known $\theta$. From this table, we see that the $L_1$, $L_2$ distances and MSE's of the WKE, TBE and Jones estimates decrease as $n$ increases while those of the naive estimate do not. Also, the $L_1$, $L_2$ distances and MSE's of the WKE and TBE are much smaller than those of the Jones estimate, which is consistent to the findings in the previous subsection that the WKE and TBE of $f$ based on estimated $\theta$ perform better than the ones based on known $\theta$.

TABLE 1
*Comparisons of estimates*

| Population | | $L_1$ | | $L_2$ | | MSE | |
| size | | mean | sd | mean | sd | mean | sd |
|---|---|---|---|---|---|---|---|
| | tbe | 0.279 | 0.126 | 0.066 | 0.056 | 0.037 | 0.034 |
| 50 | wke | 0.277 | 0.135 | 0.064 | 0.056 | 0.037 | 0.034 |
| | jones | 0.345 | 0.133 | 0.091 | 0.064 | 0.051 | 0.039 |
| | naive | 4.837 | 1.677 | 1.999 | 0.900 | 1.080 | 0.422 |
| | tbe | 0.239 | 0.100 | 0.045 | 0.037 | 0.024 | 0.022 |
| 100 | wke | 0.260 | 0.126 | 0.051 | 0.045 | 0.028 | 0.025 |
| | jones | 0.345 | 0.109 | 0.079 | 0.046 | 0.041 | 0.026 |
| | naive | 3.023 | 0.868 | 1.340 | 0.495 | 0.651 | 0.215 |
| | tbe | 0.192 | 0.072 | 0.029 | 0.021 | 0.014 | 0.012 |
| 200 | wke | 0.220 | 0.105 | 0.035 | 0.030 | 0.019 | 0.017 |
| | jones | 0.334 | 0.083 | 0.068 | 0.031 | 0.035 | 0.017 |
| | naive | 1.808 | 0.430 | 0.882 | 0.272 | 0.406 | 0.111 |
| | tbe | 0.156 | 0.045 | 0.019 | 0.010 | 0.009 | 0.006 |
| 400 | wke | 0.179 | 0.091 | 0.024 | 0.022 | 0.013 | 0.012 |
| | jones | 0.328 | 0.060 | 0.063 | 0.021 | 0.032 | 0.012 |
| | naive | 0.613 | 0.178 | 0.261 | 0.110 | 0.170 | 0.077 |
| | tbe | 0.140 | 0.039 | 0.016 | 0.008 | 0.007 | 0.004 |
| 800 | wke | 0.148 | 0.073 | 0.016 | 0.016 | 0.009 | 0.008 |
| | jones | 0.322 | 0.044 | 0.058 | 0.015 | 0.030 | 0.008 |
| | naive | 1.109 | 0.207 | 0.647 | 0.176 | 0.185 | 0.023 |
| | tbe | 0.124 | 0.033 | 0.012 | 0.006 | 0.005 | 0.003 |
| 1600 | wke | 0.126 | 0.060 | 0.012 | 0.011 | 0.006 | 0.006 |
| | jones | 0.323 | 0.032 | 0.058 | 0.011 | 0.030 | 0.006 |
| | naive | 0.838 | 0.317 | 0.396 | 0.218 | 0.191 | 0.028 |
| | tbe | 0.118 | 0.034 | 0.011 | 0.006 | 0.005 | 0.003 |
| 3000 | wke | 0.107 | 0.048 | 0.009 | 0.008 | 0.004 | 0.004 |
| | jones | 0.319 | 0.024 | 0.056 | 0.008 | 0.028 | 0.005 |
| | naive | 0.855 | 0.263 | 0.425 | 0.177 | 0.199 | 0.021 |

**95% Confidence Bands (WKE)**

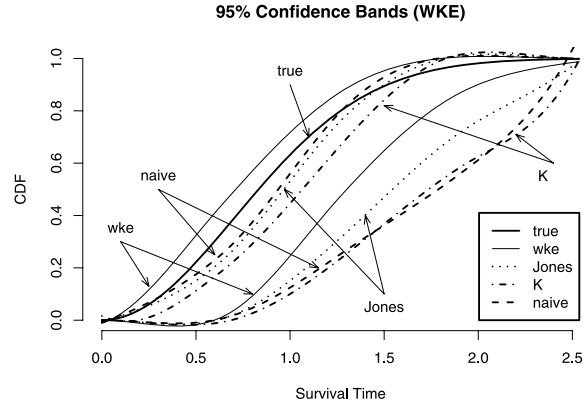

FIG 5. *Weighted kernel estimates of f.*

**95% Confidence Bands (TBE)**

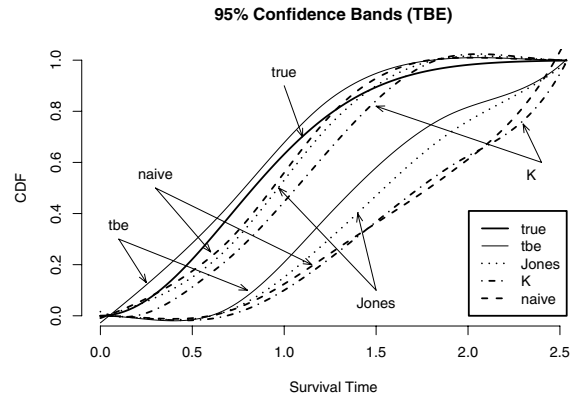

FIG 6. *Weighted kernel estimates of f.*

In our second experiment, we take $\theta = 0.5$ and repeat step 1 through step 3 and then compute the 95% pointwise confidence bands (based on the 2.5 and 97.5 percentage points from the repeats for each point) for the TBE, WKE, Jones estimate, the Kaplan–Meier estimate and the naive estimate. In this case, the length-biased assumption ($\theta = 1$) assumed in the Jones estimate is off from the true $\theta = 0.5$. From Figure 5 and Figure 6, we can easily find that when both selection bias and censoring bias exist, only our 95% confidence bands from TBE and WKE cover completely the true CDF of the survival times (the solid curve in the middle). The Jones estimate, the Kaplan–Meier estimate, and the naive estimate under-estimated $F(t)$ substantially. In Figure 5, we select the constant $c$ by the criteria in (14) and in Figure 6, we select the constant $c$ by the criteria in (13).

### 4.4. *Remarks*

A direct plug-in methodology was used to select the bandwidth in our study as that used in [15, 16, 23]. This algorithm was built into an **R** package: *KernSmooth 2.22*

For simplicity, we used the Gaussian kernel here. Other kernels such as the boundary kernels can also be used to correct the boundary effect for the survival data (survival times will never be negative.)

Some rough knowledge about $\kappa$ or $\theta$ might be used to restrict the range of the search for $c$. Here we restrict $c \in (1, 20)$.

## 5. An application

In a data set of bone marrow transplantation, there are a total of 137 patients, who were treated at one of four hospitals. The study involves transplants conducted at these institutions from March 1, 1984, to June 30, 1989. The maximum follow-up was 7 years. There were 42 patients who relapsed and 41 who died while in remission. Twenty-six patients had an episode of acute Graft-versus-host disease, and 17 patients either relapsed or died in remission without their platelets returning to normal levels [10].

Several potential risk factors were measured at the time of transplantation. For each disease, patients were grouped into risk categories based on their status at the time of transplantation. The categories were as follows: acute lymphoblastic leukemia (*ALL*), acute myelocytic leukemia (*AML*) low-risk and *AML* high-risk. Here we will focus on the disease-free survival probabilities for *ALL*, *AML* low-risk and *AML* high-risk patients. An individual is said to be disease-free at a given time after transplant if that individual is alive without the recurrence of leukemia.

There are 38 patients in group *ALL*, 54 patients in group *AML* low-risk and 45 patient in group *AML* high-risk. Figure 7 shows the estimates of the cumulative distribution function $F(t)$ with the K–M estimator for the three groups. Because the largest times in the study are different for these three groups, we find that these three estimates end at different points. The figure also suggests that if no sampling bias exists in the data, the patients in group *AML* low-risk have the most favorable prognosis (dashed line), while the *AML* high-risk group has the least favorable prognosis (dash-dotted line).

Now we use the new method to estimate from the above data by considering the possible selection bias. Here, six different estimates of $\theta$ can be obtained by TBE and WKE and by using the three different cross-validation methods in (12)–(14). Sometimes, we will not get exactly the same value by these six different methods. Which one shall we use? Simulation studies were performed for this purpose. We first generate data sets similar to the data of the three groups – with similar distributions and same proportion of data points censored. Second, these three data sets are resampled under selection functions $w'(t, \theta)$ with different $\theta$ values ($\theta = 0.5, 1, 1.5$) and the sample sizes are taken to be identical to the sample sizes of the data of groups ALL, AML low-risk and AML high-risk. Finally, $\theta$'s were estimated from those simulated samples with the new estimators (six different combinations). The above procedure was repeated 100 times for each fixed $\theta$ value. We found that the WKE with the CV$_3$ criterion defined in (14) was the winner for all $\theta$ values. For different applications, the conclusions may vary. The estimated selection functions for these three groups are

$$
\begin{array}{lll}
(16) & w(t) \propto t^{0.45} & ALL; \\
(17) & w(t) \propto t^{0.89} & AML\ low\ risk; \\
(18) & w(t) \propto t^{0.89} & AML\ high\ risk.
\end{array}
$$

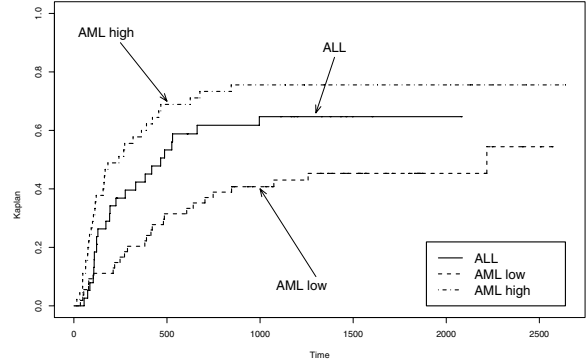

FIG 7. *K–M estimate of CDF of survival times.*

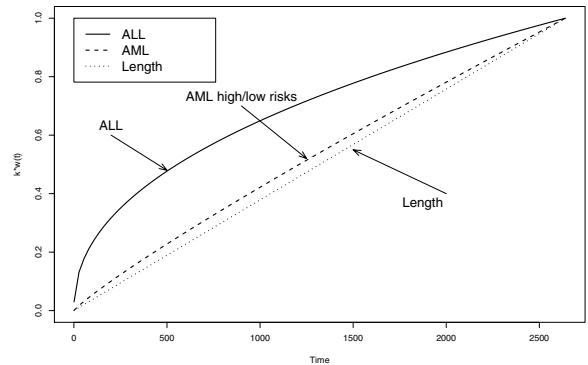

FIG 8. *Selection functions for the bone marrow transplantation data.*

These functions are plotted in Figure 8. From Figure 8, we can see that the biased sampling scheme for group *ALL* is different from those of the other two groups. In groups *AML lower-risk* and *AML high-risk*, patients with longer survival times without the recurrence of leukemia are relatively more likely to be involved in studies, $\theta = 0.89$. The selection functions are close to the selection function of the length biasing (dotted straight line, in which $\theta = 1$). While the group *ALL* has a relatively flatter selection function with larger survival times ($\theta = 0.45$). Without considering the selection bias, the actual risks will be under-estimated in all three groups.

By considering the effects caused by the biased sampling, the new cumulative failure probabilities for patients in the three groups are computed and are shown in Figure 9. From Figure 9, we can find that the risks of the patients in group *AML high-risk* are higher than those of the other two groups. This is consistent with the result of Kaplan–Meier estimates. What differs from the K–M estimate is that the risk of group *AML low-risk* is actually as high as that of group *ALL* at least in the early stage.

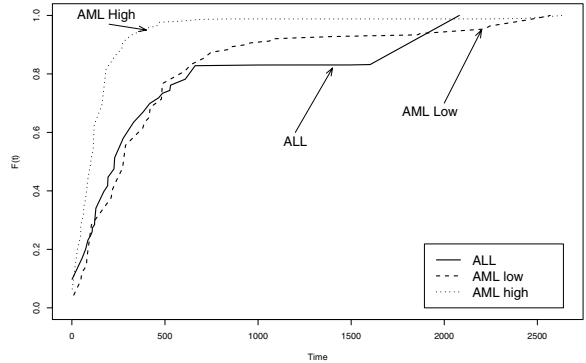

Fig 9. *New results for $\widehat{F}(t)$.*

## 6. Discussion

Since our estimation procedure also allows for constant selection function, our semi-parametric model for biased survival data is more general than the standard model. So our estimates can be used as a general-purpose procedure for analyzing survival data with a right censoring scheme that may or may not be subject to biased sampling. If our estimates are much different from the standard estimates such as the K–M estimate that ignores the selection bias, then the MAR assumption or no biased sampling assumption may be invalid, and caution must be exercised in using the standard estimate for biased survival data, which could be misleading.

In our simulation experiments, we considered $w(t) \propto t^{\theta}$ for a one-dimensional $\theta$. The resulting semi-parametric model is more general than the length-biased or area-biased sampling models. The procedure should also work for other parametric forms of $w$ and/or for multidimensional $\theta$ as long as Model Assumption 2.1 is valid. In practice, which family of $w$ should we use? Some empirical experience may help us in choosing such a family. Research on model selection of $w$ is needed. In the absent of either of these two aids, we recommend to start from a polynomial family of $w$ for some reasonable range of $\theta$.

We used a kernel density estimate to estimate $g$ in the TBE or $f_w$ in the WKE, and the Kaplan–Meier estimate to account for censoring bias. Other nonparametric smoothing estimates of density, and nonparametric estimates other than the Kaplan–Meier estimate of the survival function can in principle also be used in building our new estimates for $f$ and $\theta$.

A full-fledged asymptotic analysis of our estimators is fairly difficult and is not the objective of this paper. However, heuristically, if $\theta$ is known, it is conceivable that the TBE and WKE are consistent to $f$. When $\theta$ is unknown, if the sieve likelihood in (9) is a smooth function of $\theta$, then the plug-in estimate of $f$ by a good estimate $\hat{\theta}$ can be shown to be consistent. Note that we do not really need $\hat{\theta}$ to be consistent; all we need is that the resulting estimated selection function $\hat{w}(t, \hat{\theta})$ is consistent to $w(t, \theta)$ up to a proportional constant at $t = t_i$'s. See (1) and the expressions of the WKE and the TBE. We conjecture that the plug-in estimate of $f$ by the pseudo mle of $\theta$ is consistent under the sieve identifiability condition. This conjecture is supported by the general asymptotic property of sieve estimates (see, e.g. Bickel et al. [1]) and is confirmed by simulation results shown in Table 1.

**Acknowledgment**

Thanks to three referees and an AE for their valuable comments that helped greatly to improve this paper.